\RequirePackage{ifpdf}
\ifpdf 
\documentclass[pdftex]{sigma}
\else
\documentclass{sigma}
\fi

\newfont{\bbd}{msbm10 scaled\magstep1}
\def\Zbbd{\hbox{\bbd Z}}
\def\a{\alpha}
\def\b{\beta}
\def\g{\gamma}
\def\G{\Gamma}
\def\d{\delta}

\def\l{\lambda}

\def\pn{\phantom}
\def\ds{\displaystyle}\def\blacksquare{\rule{2mm}{2mm}}
 
\def\p{\partial}

\def\id{\hbox{{1}\kern-.25em\hbox{\rm l}}}
\def\one#1{#1^{\raise5pt\hbox{$\scriptstyle\!\!\!\!1$}}\,{}}
\def\two#1{#1^{\raise5pt\hbox{$\scriptstyle\!\!\!\!2$}}\,{}}

\def\comment#1{}

\def\?{(?)\marginpar{|?}}

\def\beq{\begin{equation}}
\def\eeq{\end{equation}}
\def\be{\begin{displaymath}}
\def\ee{\end{displaymath}}
\def\bea{\begin{eqnarray}}
\def\eea{\end{eqnarray}}
\def\bmat{\left(\begin{array}}
\def\emat{\end{array}\right)}
\newtheorem{theo}{Theorem}

\newcounter{subequation}[equation]
\makeatletter \expandafter\let\expandafter \reset@font\csname
reset@font\endcsname
\def\subeqnarray{\arraycolsep1pt
    \def\@eqnnum\stepcounter##1{\stepcounter{subequation}%
        {\reset@font\rm(\theequation\alph{subequation})}}
\jot5mm     \eqnarray}

\makeatother

\begin{document}

\numberwithin{equation}{section}

\allowdisplaybreaks

\renewcommand{\PaperNumber}{014}

\FirstPageHeading

\renewcommand{\thefootnote}{$\star$}

\ShortArticleName{An Analytic Formula for the $A_2$ Jack
Polynomials}

\ArticleName{An Analytic Formula for the $\boldsymbol{A_2}$ Jack
Polynomials\footnote{This paper is a contribution to the Vadim
Kuznetsov Memorial Issue ``Integrable Systems and Related
Topics''. The full collection is available at
\href{http://www.emis.de/journals/SIGMA/kuznetsov.html}{http://www.emis.de/journals/SIGMA/kuznetsov.html}}}

\Author{Vladimir V.~MANGAZEEV}

\AuthorNameForHeading{V.V.~Mangazeev}

  \Address{Department of Theoretical Physics,
  Research School of Physical Sciences and Engineering, \\
  The Australian National University,
  Canberra, Australia}
  \Email{\href{mailto:vladimir@maths.anu.edu.au}{vladimir@maths.anu.edu.au}}

\ArticleDates{Received November 01, 2006, in f\/inal form January
05, 2007; Published online January 24, 2007}

\Abstract{In this letter I shall review my joint results with
Vadim Kuznetsov and Evgeny Sklyanin [{\it Indag. Math.} {\bf 14}
(2003), 451--482] on separation of variables (SoV) for the $A_n$
Jack polynomials. This approach originated from the work [{\it
RIMS Kokyuroku} {\bf 919} (1995), 27--34] where the integral
representations for the $A_2$ Jack polynomials was derived. Using
special polynomial bases I shall obtain a more explicit expression
for the $A_2$ Jack polynomials in terms of generalised
hypergeometric functions.}


\Classification{05E05; 33C20; 82B23}

\vspace{-2mm}

\begin{flushright}
\it Dedicated to the memory of Vadim Kuznetsov
\end{flushright}
\vspace{-3.5mm}

\section{Introduction}

\vspace{-1mm}

In middle 1990's Vadim Kuznetsov and Evgeny Sklyanin started to
work on applications of the Separation of Variables (SoV) method
to symmetric polynomials. When I read their papers on SoV for Jack
and Macdonald polynomials, I was quite charmed by a beauty of this
approach. It connected together few dif\/ferent areas of
mathematics and mathematical physics including symmetric
polynomials,  special functions and  integrable models.

Later on we started to work together on development of SoV for
symmetric functions. During my visits to Leeds around 2001--2003
we realised that the separating operator for the $A_n$ Jack
polynomials can be constructed using Baxter's $Q$-operator
following a general approach \cite{KS5} developed by Vadim and
Evgeny Sklyanin.

Vadim was a very enthusiastic and active person. His creative
energy always encoura\-ged me to work harder on joint projects.
His area of expertise ranged from algebraic geometry to condensed
matter physics. He was an excellent mathematician and a very
strong physi\-cist.

During my visits to Leeds University I was charmed by the English
countryside. I will never forget long walks in the forest with
Vadim's family and quiet dinners at his home.

The sad news about his death was a real shock for me. His untimely
death is a big loss for the scientif\/ic community. However, he
will be long remembered by many people for his outstanding
scientif\/ic achievements.

\section[Quantum Calogero-Sutherland model]{Quantum Calogero--Sutherland model}

The quantum Calogero--Sutherland model \cite{Cal,Suth} describes a
system of $n$ particle on a circle
 with coordinates $0\le q_i\le \pi$, $i=1,\ldots,n$.
 The Hamiltonian and momentum are
 \beq
 H=-\frac{1}{2}\sum_{i=1}^n\frac{\p^2}{\p q_i^2}+\sum_{i<j}
 \frac{g(g-1)}{\sin^2(q_i-q_j)},\qquad P=-i\sum_{j=1}^n\frac{\p}{\p q_j}.
\label{jack1}
 \eeq
The ground state of the model is given by \beq
\Omega(\boldsymbol{q})=\bigg[\prod_{i<j}\sin(q_i-q_j)\bigg]^g,\qquad
\boldsymbol{q}\equiv(q_1,\ldots,q_n)\label{jack2}
 \eeq
with the ground state energy $E_0=\frac{1}{6}g^2(n^3-n)$. Further
we shall assume that $g>0$ that simplif\/ies description of the
eigenvectors.

The Calogero--Sutherland model is completely integrable \cite{OP}
and there is a commutative ring of dif\/ferential operators $H_i$,
$i=1,\ldots,n$ \cite{Sek}, which contains the Hamiltonian and
momentum~(\ref{jack1}).

The eigenfunctions are labelled by partitions
$\l\equiv(\l_1,\ldots,\l_n)$, $\l_i\ge\l_{i+1}\ge0$,
$i=1,\ldots,n-1$  and can be written as \bea
\Psi_\l(\boldsymbol{q})=\Omega(\boldsymbol{q})P_\l^{(1/g)}(\boldsymbol{x}),\label{jack3}
\eea where $P_\l^{(1/g)}(\boldsymbol{x})$,
$\boldsymbol{x}\equiv(x_1,\ldots,x_n)$ are symmetric Jack
polynomials in $n$ variables \cite{Jack,Macd,Stan} \beq
x_i=e^{2iq_i},\qquad i=1, \ldots,n. \label{jack4} \eeq

Further for a partition $\l=(\l_1,\ldots,\l_n)$,
 we def\/ine its length
$l(\l)=n$, its weight $|\l|=\l_1+\cdots+\l_n$ (see \cite{Macd})
and use a short notation $\l_{ij}=\l_i-\l_j$.

Conjugating the Hamiltonian $H$ with vacuum (\ref{jack2}) we
obtain
\begin{gather}
\Omega^{-1}(\boldsymbol{q})H\Omega(\boldsymbol{q})=\frac{1}{2}H_g+E_0,\label{jack5}
\\
H_g=\sum_{i=1}^n\left(x_i\frac{\partial}{\partial x_i}\right)^2+
g\sum_{i<j}\frac{x_i+x_j}{x_i-x_j}\left(x_i\frac{\partial}{\partial
x_i}- x_j\frac{\partial}{\partial x_j}\right).\label{jack6}
\end{gather}

Then we obtain  that Jack polynomials
$P_\l^{(1/g)}(\boldsymbol{x})$ are
 eigenfunctions of the operator $H_g$ \beq H_gP_\l^{(1/g)}(\boldsymbol{x})=E_gP_\l^{(1/g)}(\boldsymbol{x})\label{jack7} \eeq and the eigenvalues
are given by \beq E_g=\sum_{i=1}^n
\l_i[\l_i+g(n+1-2i)].\label{jack8} \eeq

It follows that symmetric Jack polynomials are orthogonal with
respect to  the scalar product
\begin{gather}
(P_\l,P_\mu)=\frac{1}{(2\pi
i)^n}\oint\limits_{|x_1|=1}\frac{dx_1}{x_1}
\cdots\oint\limits_{|x_n|=1}\frac{dx_n}{x_n}\left\{\prod_{i\neq j}
(1-x_i/x_j)\right\}^g\nonumber\\
\phantom{(P_\l,P_\mu)=}{} \times{\overline{P_\l^{(1/g)}(x)}}
P_\mu^{(1/g)}(x)=0\qquad \mbox{if}\quad \l\neq\mu,\label{jac9}
\end{gather}
where $\overline x$ stands for the complex conjugate of $x$ and
due to conditions $g>0$ and (\ref{jack4}), where $q_i$~are real,
we have ${\overline{P_\l^{(1/g)}(x)}}=P_\l^{(1/g)}(x^{-1})$.

Jack polynomials are homogeneous of the degree $|\l|$ \beq
P_\l^{(1/g)}(x_1,\ldots,x_n)=x_n^{|\l|}P_\l^{(1/g)}(x_1/x_n,
\ldots,x_{n-1}/x_n,1). \label{jack10} \eeq In fact, (\ref{jack10})
corresponds to a centre-of-mass separation in the
 Calogero--Sutherland model.

Under simultaneous shift of all parts of the partition $\l$ by an
integer the Jack polynomials undergo a simple multiplicative
transformation, which can be written as \beq
P^{(1/g)}_{\l_1,\ldots,\l_n}=(x_1\ldots x_n)^{\l_n}
P^{(1/g)}_{\l_1-\l_n,\ldots,\l_{n-1}-\l_n,0}(x_1,\ldots,x_n).\label{jack11}
\eeq

We shall use properties (\ref{jack10}), (\ref{jack11}) in our
construction of the $A_2$ Jack polynomials.

\section{Separation of variables for Jack polynomials}

In this section we shall brief\/ly review basic facts about
separation of variables \cite{Skl32,Skl38} for Jack polynomials
\cite{KS4,KMS}. Our main reference is \cite{KMS}.

We shall start with def\/inition of separated polynomials. For
each partition $\l$ def\/ines a~function $f_\l(x)$ \cite{KS4,KS99}
as
\begin{equation}
f_\l(y)=y^{\l_n}(1-y)^{1-ng}
\phantom{|}_nF_{n-1}\left(\begin{array}{l}
a_1,\ldots,a_n;\\
{b_1,\ldots,b_{n-1}}\end{array}\Bigr|\>\>y\right) ,\label{sep1}
 \end{equation}
where the hypergeometric function $\pn{|}_nF_{n-1}$
\cite{Bateman,Sl} is def\/ined as \beq
\phantom{|}_nF_{n-1}\left(\begin{array}{l}
a_1,\ldots,a_n;\\
{b_1,\ldots,b_{n-1}}\end{array}\Bigr|\>\>y\right)
=\sum_{m=0}^\infty
\frac{(a_1)_m\cdots(a_n)_m}{(b_1)_m\cdots(b_{n-1})_m}\frac{y^m}{m!},\qquad
|y|<1,\label{sep2}
\end{equation}
the Pocchammer symbol $(a)_n=\frac{\G(a+n)}{\G(a)}$ and parameters
$a_i$, $b_i$ are related to the partition
$\l=\{\l_1,\ldots,\l_n\}$ as follows
\begin{equation}
a_i=\l_n-\l_i+1-(n-i+1)g,\qquad b_j=a_j+g. \label{sep3}
\end{equation}
The important property of this function is described by the
following
\begin{theo}\label{theorem1}
$f_\l(y)$ is a polynomial in $y$ of the cumulative degree $\l_1$
and also has another representation
\begin{equation}
f_\l(y)=b_\l \>y^{\l_n}\sum_{k_1=0}^{\l_{1,2}}
\cdots\sum_{k_{n-1}=0}^{\l_{n-1,n}}
\prod_{i=1}^{n-1}(1-y)^{k_i}\frac{(-\l_{i,i+1})_{k_i}}{k_i!}
\frac{(ig)_{k_1+\cdots+k_i}}{((i+1)g)_{k_1+\cdots+k_i}},
\label{sep4}
\end{equation}
where
\begin{equation}
b_\l=\prod_{i=1}^{n-1}
\frac{((n-i+1)g)_{\l_{i,n}}}{((n-i)g)_{\l_{i,n}}}.\label{sep5}
\end{equation}
\end{theo}
This theorem proved in  \cite{KMS} allows to evaluate
$f_\l(1)=b_\l$ and we shall use this fact to calculate a correct
normalisation of the $A_2$ Jack polynomials.

In fact, we can directly expand (\ref{sep1}) into the product of
two series in $y$ and multiply them to produce an expansion of the
form \beq f_\l(y)=\sum_{k=\l_n}^{\l_1}y^k \xi_k(\l;g)\label{sep6}
\eeq with coef\/f\/icients $\xi_k(\l;g)$ being expressed in terms
of ${}_{n+1}F_n$ hypergeometric functions (see \cite{KS99} for the
case of Macdonald polynomials). However, the coef\/f\/icients of
$y^k$ will be, in general, inf\/inite series depending on $\l$ and
$g$. The Theorem~\ref{theorem1} shows that, in fact, all those
coef\/f\/icients are rational functions and the expansion really
truncates at $y^{\l_1}$. It is not clear how to produce an
explicit form of those coef\/f\/icients in terms of hypergeometric
functions of one variable with f\/inite number of terms. It can be
easily done for the case $n=2$ and later we will show how to do
that for the case $n=3$.

Now we shall formulate two main theorems of \cite{KMS}.
\begin{theo}\label{theorem2}
There exists an integral operator ${\mathcal S}_n$, which maps
symmetric polynomials in $n$~va\-riab\-les into symmetric
polynomials, \beq {\mathcal S}_n{\bf 1}={\bf 1}\label{sep7} \eeq
and this operator factorises Jack polynomials \beq {\mathcal
S}_n[P_{\l}^{(1/g)}]({x_1,\ldots,x_n})=c_\l b_\l^{-n}\prod_{i=1}^n
f_\l(x_i), \label{sep8} \eeq where \beq
c_\l=P_\l^{(1/g)}(1,\ldots,1)=\prod_{1\le i<j\le
n}\frac{(g(j-i+1))_{\l_{i,j}}} {(g(j-i))_{\l_{i,j}}}.\label{sep9}
\eeq
\end{theo}
Action  of this operator ${\mathcal S}_n$ on elementary symmetric
polynomials $e_k$ was described in~\cite{KMS}. This symmetric
separation of variables corresponds to the so called dynamical
norma\-li\-zation of the Baker--Akhiezer function (see, for
example, \cite{KS99a}).

However, using the property (\ref{jack10}) we can trivially
separate one variable $x_n$. After that we need to separate the
$n-1$ remaining variables. It results in the so called standard
separation
\begin{theo}\label{theorem3}
There exist an integral operator ${\mathcal S}^s_n$ which maps
symmetric polynomials in $n$ va\-riab\-les into polynomials
symmetric only in $n-1$ variables, \beq {\mathcal S}^s_n{\bf
1}={\bf 1}\label{sep10} \eeq and \beq {\mathcal
S}^s_n[P_{\l}^{(1/g)}]({x_1,\ldots,x_n})=c_\l b_\l^{1-n}
|x_n|^\l\prod_{i=1}^{n-1} f_\l(x_i). \label{sep11} \eeq
\end{theo}

In two next sections we shall describe the action of ${\mathcal
S}_2$ and ${\mathcal S}_3^s$ explicitly.

\section[The $A_1$ case]{The $\boldsymbol{A_1}$ case}

In this section we shall consider the simplest $n=2$ case in
standard and dynamical normalisations.

In the standard normalisation the separating operator ${\mathcal
S}_2^s$ is trivial, $c_\l=b_\l$, $\l=(\l_1,\l_2)$ and we simply
obtain \beq
P_\l^{(1/g)}(x_1,x_2)=x_2^{\l_1+\l_2}f_{\l}(x_1/x_2),\label{a1}
\eeq where the separating polynomial
\begin{gather}
f_{\l}(x)=x^{\l_2}(1-x)^{1-2g}\phantom{|}_2F_1(-\l_{12}+1-2g,1-g;
-\l_{12}+1-g;x)\nonumber\\
\phantom{f_{\l}(x)}{}
=x^{\l_2}\phantom{|}_2F_1(-\l_{12},g;-\l_{12}+1-g;x)\label{a2}
\end{gather}
and we used a well known formula for ${}_2F_1$ \cite{Bateman}. The
last formula in (\ref{a2}) gives an explicit expansion of
$f_\l(x)$ in series of $x$ with rational coef\/f\/icients.

Although formulas (\ref{a1}), (\ref{a2}) give explicit
representation for the $A_1$ Jack polynomials, a~symmetry w.r.t.\
$x_1$ and $x_2$ is not obvious. To obtain such symmetric
representation it is instructive to
 construct a separation  which corresponds to the dynamical normalisation
 of the Baker--Akhiezer function.

Introduce the following basis in the space of symmetric
polynomials $S(x_1,x_2)$ in two variab\-les~$x_1$, $x_2$ \beq
p_{mn}\equiv p_{mn}(x_1,x_2)=(x_1x_2)^m[(1-x_1)(1-x_2)]^n,\qquad
m,n\ge0  \label{a3} \eeq and def\/ine the  operator ${\mathcal
S}_2$ acting on $S(x_1,x_2)$ as \beq {\mathcal
S}_2[p_{mn}]=\frac{(g)_n}{(2g)_n}p_{mn}.\label{a4} \eeq Such
operator ${\mathcal S}_2$ can be easily represented as
 an integral operator
with a simple kernel (see \cite{KS99} for Macdonald polynomials)
but we do not need here its explicit form. What is more important
is that it maps the $A_1$ Jack polynomials into the product of
separated polynomials \beq {\mathcal
S}_2[P_{\l}^{(1/g)}](x_1,x_2)=\frac{(g)_{\l_{12}}}
{(2g)_{\l_{12}}}f_{\l}(x_1)f_{\l}(x_2).\label{a5} \eeq The action
of the inverse operator ${\mathcal S}_2^{-1}$ on the basis
$p_{mn}$ is easily constructed from (\ref{a4}) and to obtain a
symmetric representation of $P_{\l}^{(1/g)}(x_1,x_2)$ we have to
expand the rhs of (\ref{a5}) in $p_{mn}$.

To do that we use Watson's formula \cite{Wa}
\begin{gather}
\phantom{|}_2F_1(-n,b;c;x_1)\phantom{|}_2F_1(-n,b;c;x_2)\nonumber\\
\qquad
{}=\frac{(c-b)_n}{(c)_n}F_4(-n,b;c,1-n+b-c;x_1x_2,(1-x_1)(1-x_2)),\qquad
n\in\Zbbd_+ \label{a6}
\end{gather}
and the Appell's function $F_4$ is def\/ined as \beq
F_4(a,b;c,d;x,y)=\sum_{m,n=0}^\infty\frac{(a)_{m+n}(b)_{m+n}}{(c)_m(d)_n}
\frac{x^my^n}{m!n!}.\label{a7} \eeq After little algebra we obtain
from (\ref{a2})--(\ref{a7}) \beq P_{\l}^{(1/g)}(x_1,x_2)=
(x_1x_2)^{\l_2}\sum_{0\le m+n\le\l_{12}}
\frac{(-\l_{12})_{m+n}(g)_{m+n}}{(-\l_{12}+1-g)_m(g)_n}
\frac{p_{mn}(x_1,x_2)}{m!n!}.\label{a8} \eeq Finally we can
rewrite (\ref{a8}) in the basis of elementary symmetric
polynomials $e_1=x_1+x_2$, $e_2=x_1x_2$. Expanding $p_{mn}$ in
$e_1$, $e_2$ and evaluating remaining sums we come to \beq
P_{\l}^{(1/g)}(x_1,x_2)= (e_2)^{\l_2}(-1)^{\l_{12}}
\frac{(\l_{12})!}{(g)_{\l_{12}}}
\sum_{i+2j=\l_{12}}(g)_{i+j}\frac{e_1^ie_2^j}{i!j!}(-1)^{i+j}.\label{a9}
\eeq Another symmetric representation for
$P_{\l}^{(1/g)}(x_1,x_2)$ which f\/irst appeared in \cite{Koorn}
is given in terms of Gegenbauer polynomials $C_n^g(x)$. Expanding
$e_1$ in series of $x_1$ and $x_2$ we obtain \beq
P_{\l}^{(1/g)}(x_1,x_2)=
(x_1x_2)^{\frac{1}{2}|\l|}\frac{(\l_{12})!}{(g)_{\l_{12}}}
C_{\l_{12}}^g\left(
\frac{1}{2}\big[(x_1/x_2)^{1/2}+(x_2/x_1)^{1/2}\big]\right).\label{a10}
\eeq In fact, this formula can be directly obtained from
representation (\ref{a1}), (\ref{a2}), but we preferred to use the
symmetric separation (\ref{a5}). As was shown in \cite{KS99a}, the
separating operator ${\mathcal S}_{n-1}$ for the dynamical
normalisation has the same structure as the separating operator
${\mathcal S}_n^s$ in the standard normalisation. In the next
section we shall use this fact to construct a representation for
the $A_2$ Jack polynomials.

\section[The $A_2$ case and the expansion of products of
separated polynomials]{The $\boldsymbol{A_2}$ case and the
expansion of products\\ of separated polynomials} 

Now we shall consider $n=3$ case in the standard normalisation.
Using (\ref{jack10}) we obtain \beq
P_\l^{(1/g)}(x_1,x_2,x_3)=x_3^{|\l|}P_\l^{(1/g)}(x_1/x_3,x_2/x_3,1)\equiv
x_3^{|\l|}p_{\l}(x_1/x_3,x_2/x_3),\label{part1} \eeq where we
introduced a short notation $p_\l(x_1,x_2)$ for the $A_2$ Jack
polynomials with $x_3=1$.

Def\/ine a reduced separating operator $\hat{\mathcal S}_3$ which
acts only on the basis $p_{mn}$ (see (\ref{a3})) in the space of
symmetric polynomials in two variables $S(x_1,x_2)$ as \beq
\hat{\mathcal
S}_3[p_{mn}]=\frac{(2g)_n}{(3g)_n}p_{mn}.\label{part2} \eeq The
explicit construction of $\hat{\mathcal S}_3$ as an integral
operator was given in \cite{KS4,KMS}. The operator $\hat{\mathcal
S}_3$ separates variables in $p_\l(x_1,x_2)$ \beq \hat{\mathcal
S}_3[p_\l](x_1,x_2)= c_\l
b_\l^{-2}f_\l(x_1)f_\l(x_2),\label{part3} \eeq where $f_\l(x)$ is
given by (\ref{sep1}) with $n=3$.

To construct $p_\l(x_1,x_2)$ we have to invert the action of
$\hat{\mathcal S}_3$. It follows from (\ref{part2}) that in order
to apply to (\ref{part3}) $\hat{\mathcal S}_3^{-1}$ we have to
expand the rhs of (\ref{part3}) in  $p_{mn}$ \beq
f_{\l}(x_1)f_{\l}(x_2)=(x_1x_2)^{\l_3} \sum_{0\le
m+n\le\l_{1,3}}c_{m,n}(\l;g)\>p_{mn}(x_1,x_2)\label{part4} \eeq
 and it is easy to see that coef\/f\/icients $c_{m,n}(\l;g)$ depend
on the partition $\l$ only via $\l_{13}$ and $\l_{23}$.

In fact, formula (\ref{part4}) is the main nontrivial result of
this paper. Setting $x_2=1$ in (\ref{part4}) we obtain as a
particular case the expansion for $f_\l(x_1)$ of the form
(\ref{sep6}) with coef\/f\/icients \beq
\xi_k(\l,g)=c_{m,\,0}(\l,g)/b_\l.\label{part4a} \eeq

Now we shall calculate coef\/f\/icients $c_{m,n}(\l,g)$. For
convenience we shall omit dependence of $c_{m,n}(\l;g)$ on $\l$
and $g$.

Let us remind that the hypergeometric function ${}_3F_{2}$
satisf\/ies $3$-th order dif\/ferential equation of the form
\cite{Bateman,Sl} \beq
\left[\d\prod_{i=1}^{2}(\d+b_i-1)-x\prod_{i=1}^{3}(\d+a_i)\right]
\phantom{|}_{3}F_{2}(a_1,a_2,a_3;b_1,b_2;x)=0, \qquad \d\equiv
x\frac{d}{dx}.\label{part5} \eeq Using (\ref{sep1})
 we can easily
obtain a third order dif\/ferential equation for the functions
$f_\l(x_i)$, $i=1,2$. Now we can rewrite these equations in terms
of new variables
 \beq u=x_1x_2,\qquad v=(1-x_1)(1-x_2),\qquad
p_{m,n}=u^mv^n.\label{part6} \eeq Therefore, for the rhs of
(\ref{part4}) we obtain two dif\/ferential equations in $u$ and
$v$ which can be rewritten as recurrence relations for the
coef\/f\/icients $c_{m,n}$. All calculations are straightforward
and we shall omit the details. The f\/irst recurrence relation for
$c_{mn}$ looks amazingly simple and basically allows to f\/ind
$c_{m,n}$ in a closed form {\samepage
\begin{gather}
 c_{m+1,n}(m+1)(1-2g+m-r_1)(1-g+m-r_2)\nonumber\\
\qquad{}{} +
c_{m,n+1}(n+1)(3g+n-1)(3g+n)\nonumber\\
\qquad{} +c_{mn}(2g+m+n)(r_1-m-n)(g+m+n-r_2)=0,\label{part7}
\end{gather}
where we denoted $r_1=\l_{1,3}$, $r_2=\l_{2,3}$.}

The relation (\ref{part7}) immediately suggests to make a
substitution \beq
c_{mn}=\frac{(2g)_{m+n}(-r_1)_{m+n}(g-r_2)_{m+n}}
{m!n!(3g-1)_n(3g)_n(1-2g-r_1)_m(1-g-r_2)_m}\>a_{mn}.\label{part8}
\eeq Then (\ref{part7}) reduces to \beq
a_{m+1,n}+a_{m,n+1}=a_{mn}\label{part9} \eeq which can be  solved
in terms of $a_{m,0}$ \beq a_{m,n}=\sum_{l=0}^n
(-1)^l\binom{n}{l}a_{m+l,0}.\label{part10} \eeq

The second recurrence relation for $c_{m,n}$ leads to the
following relation for $a_{m,n}$
\begin{gather}
 m(2g+r_1-m)(g+r_2-m)a_{m-1,n}-n(3g+n-2)(3g+n-1)a_{m,n-1}
\nonumber\\
\qquad{}+(m+n-r_1)(2g+m+n)(g-r_2+m+n)a_{m+1,n}\nonumber\\
\qquad{} +[n(3g+n-1)(3g+n-2)-3m(m+1)n-n(r_1-1)(r_2-1) \nonumber\\
\qquad{}+2(m-r_2)[g(1+r_1)+(r_1-m)m]+2(3g+r_1+r_2)mn\nonumber\\
\qquad{}-g(g-1)(r_1+3r_2-5m)-g(2g-3+r_1+2r_2)n]\>a_{mn}=0.\label{part11}
\end{gather}
However, due to (\ref{part10}) we need the latter relation only
for the calculation of $a_{m,0}$ which corresponds to the case
$n=0$ in (\ref{part11}). In fact, it is much easier to obtain the
$3$-term recurrence relation for $a_{m,0}$ not from
(\ref{part11}), but from the $3$-term recursion for the
coef\/f\/icients $c_{m,0}$ itself. This latter recursion follows
from (\ref{sep6}), (\ref{part4a}) and the third order
dif\/ferential equation for $f_\l(x)$.

Now let us introduce a generating function \beq
h(x)=x^{1-\rho}\sum_{m\ge0}a_{m,0}\>x^m\label{part12} \eeq and
rewrite the recurrence relation for $a_{m,0}$ as a third order
dif\/ferential equation for $h(x)$. Then this equation can be
easily reduced to the standard hypergeometric equation
(\ref{part5}) for the function ${}_3F_2$. This equation has three
solutions  with exponents \beq \rho=2g,g-r_2,-r_1\label{part13}
\eeq
 at $x=0$ which can be written as
\beq h(x)=x^{1-\rho}(1-x)^{3g-2}
\phantom{|}_4F_3\left(\begin{array}{r}
g-r_1-\rho,3g-\rho,2g-r_2-\rho,1;\\
{2g+1-\rho,1-\rho-r_1,1-\rho+g-r_2}\end{array}
\Bigr|\>x\right)\label{part14} \eeq and for any choice of $\rho$
a hypergeometric function ${}_4F_3$ reduces to ${}_3F_2$.

Now we note that for $\rho=2g$ and $\rho=g-r_2$ the hypergeometric
function in the rhs of~(\ref{part14}) truncates and we obtain from
(\ref{part12}), (\ref{part14}) the following equivalent
expressions for the coef\/f\/icients $a_{m,0}$ \beq
a_{m,0}=\a_1\frac{(1-g)_m}{(2g)_m}
\phantom{|}_4F_3\left(\begin{array}{c}
-r_2, \>g,\>-g-r_1,1-2g-m;\\
{1-g-r_2,1-2g-r_1,g-m}\end{array}\Bigr|\>1\right)\label{part16}
\eeq and \beq a_{m,0}=\a_2 \frac{(1-2g-r_2)_m}{(g-r_2)_{m}}
\phantom{|}_4F_3\left(\begin{array}{l}
r_2-r_1, \>g,\>2g+r_2,1-g+r_2-m;\\
{1-g+r_2-r_1,1+g+r_2,2g+r_2-m}\end{array}\Bigr|\>1\right).\label{part17}
\eeq

To calculate the correct normalisation coef\/f\/icients $\a_1$,
$\a_2$ for (\ref{part16}), (\ref{part17}) we take the limit
$x_1\to 0$, $x_2\to 1$ in (\ref{part4}). Then only the term
$c_{0,0}$ contributes to the rhs of (\ref{part4}) and we use
(\ref{sep5}), (\ref{part8}) to calculate $a_{0,0}$.

From the other side the hypergeometric series in
(\ref{part16}), (\ref{part17}) becomes Saalsch\"utzian ${}_3F_2$
(i.e.\ $1+a_1+a_2+a_3= b_1+b_2$) and can be evaluated using the so
called Saalsch\"utz's theorem \cite{Bateman,Sl} \beq
{}_3F_2\left(\begin{array}{l}
a,b,-n;\\
{c,1+a+b-c-n}\end{array}\right)=
\frac{(c-a)_n(c-b)_n}{(c)_n(c-a-b)_n}. \eeq Then we obtain
\begin{gather}
\a_1=\frac{(r_1-r_2)!}{r_1!}\frac{(3g)_{r_1}(2g)_{r_2}}
{(2g)_{r_1-r_2}(1-g)_{r_2}},\label{part18}
\\
\a_2=\frac{r_2!}{r_1!}
\frac{(1+g)_{r_1}(g)_{r_1-r_2}}{(1+g)_{r_2}(2g)_{r_1-r_2}}
\frac{(3g)_{r_1}(2g)_{r_2}}{(2g)_{r_1}(g)_{r_2}}.\label{part19}
\end{gather}
Now we can substitute (\ref{part16})--(\ref{part17}) into
(\ref{part10}). Expanding hypergeometric functions into series we
can interchange the order of summations and calculate the sum over
$l$ in (\ref{part10}). Substituting the result into (\ref{part8})
we f\/inally come to the following explicit expressions for
$c_{m,n}$
\begin{gather}
 c_{m,n}=\frac{(r_1-r_2)!}{r_1!}\frac{(3g)_{r_1}(2g)_{r_2}}
{(2g)_{r_1-r_2}(1-g)_{r_2}}
\frac{(-r_1)_{m+n}(g-r_2)_{m+n}(1-g)_m}
{m!n!(1-2g-r_1)_m(1-g-r_2)_m(3g)_n}\nonumber\\
\phantom{c_{m,n}=}{}\times \phantom{|}_4F_3\left(\begin{array}{c}
-r_2, \>g,\>-g-r_1,1-2g-m-n;\\
{1-g-r_2,1-2g-r_1,g-m}\end{array}\Bigr|\>1\right) \label{part20}
\end{gather}
or
\begin{gather}
 c_{m,n}=\frac{r_2!}{r_1!}\frac{(1+g)_{r_1}(g)_{r_1-r_2}}
{(1+g)_{r_2}(2g)_{r_1-r_2}}
\frac{(3g)_{r_1}(2g)_{r_2}}{(2g)_{r_1}(g)_{r_2}}
\frac{(1-2g-r_2)_m(-r_1)_{m+n}(2g)_{m+n}}
{(1-2g-r_1)_m(1-g-r_2)_m(3g)_n}\nonumber\\
\phantom{c_{m,n}=}{}\times
\frac{1}{m!n!}\phantom{|}_4F_3\left(\begin{array}{l}
r_2-r_1, \>g,\>2g+r_2,1-g+r_2-m-n;\\
{1-g+r_2-r_1,1+g+r_2,2g+r_2-m}\end{array}\Bigr|\>1\right).\label{part21}
\end{gather}

As we mentioned above the coef\/f\/icients of the expansion of the
separating polynomial $f_\l(x)$ in~$x$ are given by $c_{m,0}$,
i.e.\ expressed in terms of hypergeometric functions ${}_4F_3$. It
is quite surprising that the coef\/f\/icients $c_{m,n}$ of the
expansion (\ref{part4}) are still given in terms of ${}_4F_3$ by
formulas (\ref{part20}), (\ref{part21}) and do not involve higher
hypergeometric functions. This happens probably because the basis
$p_{mn}$ is very special.

\section[A representation for the $A_2$ Jack polynomials]{A representation for the $\boldsymbol{A_2}$ Jack polynomials}

Now we are ready to give explicit formulas for the $A_2$ Jack
polynomials. Applying the opera\-tor~$\hat{\mathcal S}_3^{-1}$ to
(\ref{part3}) and using (\ref{part1}), (\ref{part2}),
(\ref{part4}) we obtain two representations
\begin{gather}
P_{(\l_1,\l_2,\l_3)}^{(1/g)}(x_1,x_2,x_3)=
\frac{(g)_{\l_{23}}(2g)_{\l_{13}}}
{(g)_{\l_{12}}(1-g)_{\l_{23}}}\frac{(\l_{12})!}{(\l_{13})!}\>
(x_1x_2)^{\l_3}x_3^{\l_1+\l_2-\l_3}\nonumber\\
\qquad{}\times\sum_{0\le m+n\le\l_{13}}
\frac{1}{m!n!}\frac{(-\l_{13})_{m+n}(g-\l_{23})_{m+n}(1-g)_m}
{(1-2g-\l_{13})_m(1-g-\l_{23})_m(2g)_n}
\label{repr1}\\
\qquad{}\times\phantom{|}_4F_3\left(\begin{array}{c}
-\l_{23}, \>g,\>-g-\l_{13},1-2g-m-n;\\
{1-g-\l_{23},1-2g-\l_{13},g-m}\end{array}\Bigr|\>1\right)
\left[\frac{x_1x_2}{x_3^2}\right]^m
\left[\left(1-\frac{x_1}{x_3}\right)
\left(1-\frac{x_2}{x_3}\right)\right]^n\nonumber
\end{gather}
and
\begin{gather}
 P_{(\l_1,\l_2,\l_3)}^{(1/g)}(x_1,x_2,x_3)=
\frac{(g)_{\l_{13}+1}}
{(g)_{\l_{23}+1}}\frac{(\l_{23})!}{(\l_{13})!}\>
(x_1x_2)^{\l_3}x_3^{\l_1+\l_2-\l_3}
\nonumber\\
\label{repr2} \qquad{}\times\sum_{0\le
m+n\le\l_{13}}\frac{1}{m!n!}
\phantom{|}_4F_3\left(\begin{array}{c}
-\l_{12}, \>g,\>2g+\l_{23},1-g+\l_{23}-m-n;\\
{1+g+\l_{23},1-g-\l_{12},2g+\l_{23}-m}\end{array}\Bigr|\>1\right)
\\
\qquad{}\times \frac{(-\l_{13})_{m+n}(1-2g-\l_{23})_{m}(2g)_{m+n}}
{(1-2g-\l_{13})_m(1-g-\l_{23})_m(2g)_n}
\left[\frac{x_1x_2}{x_3^2}\right]^m\left[\left(1-\frac{x_1}{x_3}\right)
\left(1-\frac{x_2}{x_3}\right)\right]^n.\nonumber
\end{gather}

Formulas (\ref{repr1}), (\ref{repr2}) simplify when $\l_2=\l_3$
and $\l_1=\l_2$ accordingly.

First consider the case of one-row partitions $(\l,0,0)$. It is
convenient to set $x_3=1$. Then we obtain from (\ref{repr1}) \beq
P_{\l,0,0}^{(1/g)}(x_1,x_2,1)=\frac{(2g)_\l}{(g)_\l} \sum_{0\le
m+n\le\l} \frac{(-\l)_{m+n}(g)_{m+n}}{(1-2g-\l)_m(2g)_n}
\frac{(x_1x_2)^m[(1-x_1)(1-x_2)]^n}{m!n!}.\label{repr3} \eeq
Expanding the last term we can transform this formula to another
basis \beq P_{\l,0,0}^{(1/g)}(x_1,x_2,1)= \sum_{0\le
i+2j\le\l}\frac{(-1)^{i+j}}{i!j!}
\frac{(-\l)_{i+2j}\>(g)_{i+j}(g)_{\l-i-2j}}
{(g)_\l}(x_1+x_2)^i(x_1x_2)^j. \label{repr4} \eeq

Finally restoring dependence on $x_3$ we shall rewrite
(\ref{repr4}) in the basis of elementary symmetric functions
$e_1=x_1+x_2+x_3$, $e_2=x_1x_2+x_1x_3+x_2x_3$ and $e_3=x_1x_2x_3$.
\beq P_{\l,0,0}^{(1/g)}(x_1,x_2,x_3)=\frac{\l!}{(g)_\l}(-1)^\l
\sum_{i+2j+3k=\l} (-1)^{i+j+k}\frac{(g)_{i+j+k}}{i!j!k!}\>
e_1^ie_2^je_3^k.\label{repr5} \eeq

In fact, (\ref{repr5}) corresponds to the f\/irst fundamental
weight $\omega_1$ of $A_2$ and is a special case of the
Proposition 2.2 from \cite{Stan} for arbitrary number of variables
\beq P_{(r)}^{(1/g)}(x)=\frac{r!}{(g)_r} \sum_{|\mu|=r}e_\mu
\frac{(g)_{l(\mu)}}{\ds\prod_{i\ge1}m_i!}(-1)^{r-l(\mu)}\label{repr6}
\eeq for all partitions $\mu$ such that \beq
\mu=(1^{m_1}2^{m_2}\cdots),\qquad |\mu|=r,\qquad l(\mu)=\sum
m_i,\qquad e_\mu=e_1^{m_1}e_2^{m_2}\cdots. \eeq

Now from (\ref{repr2}) we can similarly calculate
$P_{\l,\l,0}^{(1/g)}(x_1,x_2,x_3)$ which corresponds to the second
fundamental weight $\omega_2$ of $A_2$ \beq
P_{(\l,\l,0)}^{(1/g)}(x_1,x_2,x_3)=
\frac{\l!}{(g)_\l}\sum_{|\mu|=2\l}e_\mu\frac{(g)_{\l-m_3}(-1)^{m_3}}
{\left(\l-\sum\limits_{i=1}^3
m_i\right)!\prod\limits_{i=1}^{2}m_i!} \label{repr8} \eeq for all
partitions $\mu=(1^{m_1}2^{m_2}3^{m_3})$, $|\mu|=2r$.

It seems that (\ref{repr8}) has the following generalisation to
$n$ variables \beq P_{(\l^{n-1},0)}^{(1/g)}(x_1,\ldots,x_n)=
\frac{\l!}{(g)_\l}\sum_{|\mu|=(n-1)\l}e_\mu\frac{(g)_{\l-m_n}(-1)^{m_n}}
{\left(\l-\sum\limits_{i=1}^n
m_i\right)!\prod\limits_{i=1}^{n-1}m_i!} \label{repr8a} \eeq for
all partitions $\mu=(1^{m_1}\cdots n^{m_n})$, $|\mu|=(n-1)\l$. We
are not aware of any appearances of~(\ref{repr8a}) in the
literature.

Note that formulas for the $A_2$ Jack polynomials were already
obtained earlier in \cite{Per}. They were generalised  in
\cite{Las} for partitions of length $3$ with any number of
variables and further in~\cite{LS} for the general $A_n$ case.

The approach used in \cite{Per,Las} is based on the inversion of
the Pieri formulas. For the $A_2$ case formulas of \cite{Per}  can
be written as
\begin{gather}
P_{(\l_1,\l_2,0)}=
\sum_{i=0}^{\min(\l_1-\l_2,\l_2)}\b_{\l_1,\l_2}^i
P_{(\l_1-\l_2-i,0,0)} P_{(\l_2-i,\l_2-i,0)},\label{pieri1}
\\
P_{(\l_1,\l_2,0)}=\sum_{i=0}^{\l_2}\g_{\l_1,\l_2}^i
P_{(\l_1+i,0,0)} P_{(\l_2-i,0,0)},\label{pieri2}
\end{gather}
where $\b_{\l_1,\l_2}^i$ and $\g_{\l_1,\l_2}^i$ are explicitly
given factorised coef\/f\/icients. Note also that (\ref{pieri2})
is a~specialisation of the result \cite{Jing} for partitions of
length $2$ and any number of variables.

It is interesting to compare (\ref{pieri1}), (\ref{pieri2}) with
our formulas (\ref{repr1}), (\ref{repr2}). Using the
expres\-sions~(\ref{repr5}) and (\ref{repr8}) we can rewrite
(\ref{pieri1}), (\ref{pieri2}) in terms of quintuple sum in the
basis of elementary symmetric functions. However, our formulas
(\ref{repr1}), (\ref{repr2}) contain only triple sums  in a
dif\/ferent basis.

We tried hard to obtain directly (\ref{repr1}), (\ref{repr2}) from
(\ref{pieri1}), (\ref{pieri2}) with no success. The corresponding
transformations can be reduced to multiple sum identities
involving a large number of Pochhammer symbols. It is absolutely
unclear how to prove them.

Finally, it would be interesting to rewrite formulas
(\ref{repr1}), (\ref{repr2}) in the basis $e_1$, $e_2$, $e_3$ like
(\ref{repr5}) for a general partition $\l=(\l_1,\l_2,\l_3)$.
However, it involves sums of truncated hypergeometric functions
${}_4F_3$ and we did not succeed at this stage to obtain any
compact expression for general partitions of length $3$.

\section{Conclusions} 
In this letter we applied  the separation of variables method to
construct  a formula for the $A_2$ Jack polynomials based on the
results of \cite{KS4}. We explicitly evaluated the
coef\/f\/icients of the expansion of the $A_2$ Jack polynomials in
a special basis.

It is a dif\/f\/icult problem to show that our formulas
(\ref{repr1}), (\ref{repr2}) can be reduced to
results of~\cite{Per,Las}. In fact, their structure is completely
dif\/ferent. Nevertheless, our approach can be useful to f\/ind a
simpler representation for Jack and Macdonald polynomials.

\subsection*{Acknowledgments}

I greatly benef\/ited in my understanding of the quantum SoV from
very stimulating discussions with Vadim Kuznetsov and Evgeny
Sklyanin. A warm atmosphere at Leeds University, critical and
creative talks with Vadim are things which will be greatly missed.

\pdfbookmark[1]{References}{ref}
\LastPageEnding


\begin{thebibliography}{99}

\footnotesize\itemsep=0pt


\bibitem{Cal} Calogero F.,
 Solution of a three-body problem in one dimension,
{\it J. Math. Phys.} {\bf 10} (1969), 2191--2196.

\bibitem{Bateman} Erd\'elyi A. (Editor),
Higher transcendental functions, Vol.~I, McGraw-Hill Book Company,
1953.

\bibitem{Jack} Jack H.,
 A class of symmetric functions with a parameter,
{\it Proc. Royal Soc. Edinburgh (A)} {\bf 69} (1970),  1--18.

\bibitem{Jing} Jing N.H., J\'ozef\/iak T.,
A formula for two-row Macdonald functions, {\it Duke Math. J.}
{\bf 67} (1992), 377--385.

\bibitem{Koorn} Koornwinder T.H., Sprinkhuizen-Kuyper I.G.,
 Generalized power series expansions
for a class of orthogonal polynomials in two variables, {\it SIAM
J. Math. Anal.} {\bf 9} (1978), 457--483.

\bibitem{KMS}Kuznetsov V.B., Mangazeev V.V., Sklyanin E.K.,
$Q$-operator and factorized separation chain for Jack polynomials,
{\it Indag. Math.} {\bf 14} (2003), 451--482,
\href{http://arxiv.org/abs/math.CA/0306242}{math.CA/0306242}.

\bibitem{KS4} Kuznetsov V.B., Sklyanin E.K.,
Separation of variables in the $A_2$ type Jack polynomials, {\it
RIMS Kokyuroku} {\bf 919} (1995), 27--34,
\href{http://arxiv.org/abs/solv-int/9508002}{solv-int/9508002}.

\bibitem{KS5} Kuznetsov V.B.,
Sklyanin E.K.,  On B\"acklund transformations for many-body
systems, {\it J.\ Phys.\ A: Math.\ Gen.} {\bf 31} (1998),
2241--2251,
\href{http://arxiv.org/abs/solv-int/9711010}{solv-int/9711010}.

\bibitem{KS99} Kuznetsov V.B., Sklyanin E.K., Factorisation
of Macdonald polynomials, in Symmetries and Integrability of
Dif\/ference Equations, {\it Lecture Note Series}, Vol.~255,
London Math. Society, Cambridge University Press, 1999, 370--384,
\href{http://arxiv.org/abs/q-alg/9703013}{q-alg/9703013}.

\bibitem{KS99a} Kuznetsov V.B., Sklyanin E.K., Separation of
variables and integral relations for special functions, {\it The
Ramanujan J.} {\bf 3} (1999), 5--35,
\href{http://arxiv.org/abs/q-alg/9705006}{q-alg/9705006}.

\bibitem{Las} Lassalle M.,
 Explicitation des polyn\^omes de Jack et de Macdonald en longueur trois, {\it C.R. Acad. Sci. Paris
 Ser. I Math.}
{\bf 333} (2001), 505--508.

\bibitem{LS} Lassalle M., Schlosser M., Inversion of the Pieri
formula for Macdonald polynomials, {\it Adv. Math.}, {\bf 202}
(2006), 289--325,
\href{http://arxiv.org/abs/math.CO/0402127}{math.CO/0402127}.

\bibitem{Macd}
 Macdonald I.G., Symmetric functions and Hall polynomials, 2nd ed.,
Oxford University Press, Oxford, 1995.

\bibitem{OP} Olshanetsky  M.A., Perelomov A.M.,
 Quantum integrable systems
related to Lie algebras, {\it Phys. Rep.} {\bf 94} (1983),
313--404.

\bibitem{Per} Perelomov A.M., Ragoucy E., Zaugg P.,
 Explicit solution of the quantum three-body
Calogero--Sutherland model, {\it J. Phys. A: Math. Gen.} {\bf 31}
(1998), L559--L565,
\href{http://arxiv.org/abs/hep-th/9805149}{hep-th/9805149}.

\bibitem{Sek} Oshima  T., Sekiguchi H.,
 Commuting families of dif\/ferential
operators invariant under the action of a Weyl group, {\it J.
Math. Sci. Univ. Tokyo} {\bf 2} (1995),  1--75.

\bibitem{Skl32} Sklyanin E.K.,
Quantum inverse scattering method. Selected topics, in Quantum
Group and Quantum Integrable Systems, Editor M.-L.~Ge,  {\it
Nankai
   Lectures in Mathematical Physics}, World Scientif\/ic,
   Singapore,
       1992, 63--97, \href{http://arxiv.org/abs/hep-th/9211111}{hep-th/9211111}.

\bibitem{Skl38} Sklyanin E.K., Separation of variables. New trends,
{\it Progr.\ Theoret. Phys.\ Suppl.} {\bf 118} (1995), 35--60,
\href{http://arxiv.org/abs/solv-int/9504001}{\mbox{solv-int/9504001}}.

\bibitem{Sl} Slater L.J., Generalized hypergeometric functions,
Cambridge University Press, 1966.

\bibitem{Stan} Stanley R.P.,  Some combinatorial properties of Jack
symmetric functions, {\it Adv. Math.} {\bf 77} (1989), 76--115.

\bibitem{Suth} Sutherland B.,
Quantum many-body problem in one dimension, I, II, {\it J. Math.
Phys.} {\bf 12} (1971), 246--250, 251--256.

\bibitem{Wa} Watson G.N.,
The product of two hypergeometric functions, {\it Proc. London
Math. Soc. (2)} {\bf 20} (1922), 191--195.
\end{thebibliography}
\end{document}